\documentclass[11pt]{amsart}

\usepackage{amssymb,amsmath}
\usepackage{bbm}
\usepackage{a4wide}

\newtheorem{theorem}{Theorem}

\newtheorem{lemma}{Lemma}

\theoremstyle{definition}
\newtheorem{remark}{Remark}

\newcommand{\ts}{\hspace{0.5pt}}
\newcommand{\nts}{\hspace{-0.5pt}}
\newcommand{\CC}{\mathbb{C}\ts}
\newcommand{\RR}{\mathbb{R}\ts}

\newcommand{\NN}{\mathbb{N}}

\newcommand{\cI}{\mathcal{I}}

\newcommand{\pa}{\varPhi^{}_{\! \nts A}}
\newcommand{\dd}{\,\mathrm{d}}
\newcommand{\ee}{\ts\mathrm{e}}
\newcommand{\one}{\mathbbm{1}}
 
\begin{document}

\title[Peano--Baker series]
{The Peano--Baker series}

\author{Michael Baake}
\author{Ulrike Schl\"{a}gel}

\address{Fakult\"at f\"ur Mathematik, Universit\"at Bielefeld, \newline
\hspace*{3.1ex}Postfach 100131, 33501 Bielefeld, Germany}


\begin{abstract} 
  This note reviews the Peano--Baker series and its use to solve the
  general linear system of ODEs. The account is elementary and
  self-contained, and is meant as a pedagogic introduction to this
  approach, which is well known but usually treated as a folklore
  result or as a purely formal tool. Here, a simple convergence result
  is given, and two examples illustrate that the series can be used
  explicitly as well.
\end{abstract}

\maketitle

\section{Introduction}

Consider the inhomogeneous linear initial value (or Cauchy) problem 
\begin{equation} \label{Cauchy} \dot{x} = A\ts x + b \ts , \quad
  \text{with } \, x(t^{}_0) = x^{}_0 \ts ,
\end{equation}
on some interval $I\subset \RR$ that contains $t_0$, with
time-dependent quantities $x$, $A$ and $b$. Here, $A(t)$ denotes a
family of matrices and $b(t)$ a vector-valued function, both of
dimension $d$ over $\RR$ (or over $\CC$). For simplicity, we assume
$A$ and $b$ to be continuous on $I$. When $[A(t),A(s)]=0$ for all
$t,s \in I$ (where we use the notation $[A,B]=AB - BA$), the unique
solution of \eqref{Cauchy} is
\begin{equation} \label{commuting-case}
    x(t) = \exp \Bigl( \ts\int_{t^{}_0}^{t} \! A(\tau) \dd \tau \Bigr)
       \biggl( x^{}_{0} + \int_{t^{}_0}^{t} \! \exp \Bigl( - \!
       \int_{t^{}_0}^{\tau} \! A(\sigma) \dd\sigma \Bigr)\,
       b (\tau) \dd\tau \biggr),
\end{equation}
which simplifies to
$x(t) = \exp \bigl( \ts\int_{t^{}_0}^{t} \!  A(\tau) \dd \tau \bigr)
\ts x^{}_{0} $ for the homogeneous case; see
\cite[Cor.~2.41]{KP}. Also, when $A(t)=A$ is constant, one has
$\,\exp \bigl( \ts\int_{t^{}_0}^{t} \! A(\tau) \dd \tau \bigr) =
\ee^{\ts t-t^{}_0)A}$ as usual; we refer the reader to \cite[Secs.~11
and 12]{A} for a more general exposition.

The solution of \eqref{Cauchy} is still unique when $[A(t),A(s)]\ne 0$
for $t\ne s$, but the solution is then given by the Peano--Baker
series (PBS); see \cite[Sec.~16.5]{I}, \cite[Sec.~7.5]{FDC} or
\cite[Sec.~1.3]{B} for background. This approach goes back to Peano
(1888; see \cite{P}), and was further developed by Baker (1905; see
\cite{Baker}). Baker's paper also contains a more detailed account of
the history of this approach.  The PBS is well-known in principle, but
not widely known, and appears mainly in engineering textbooks; compare
\cite{B,K,FH} for examples.  Quite frequently, it is used in a purely
formal manner (without convergence considerations), and it is also
sometimes claimed \cite{R} that it is of little practical use for the
solution. However, as one can learn from the biographical sketches of
Delone in \cite{Del}, one should neither discard analytical tools nor
believe repeated claims without evidence. In fact, when new problems
surface, many (old) tools turn out to be more useful than expected.

It is thus the purpose of this short note to summarise the elementary
properties of the Peano--Baker series and to give a simple and
self-contained account. In particular, we do not restrict ourselves to
the abstract fundamental system (or matrix solution), but discuss its
explicit series expansion with convergence. The need for it came up in
the context of a problem in mathematical population genetics; see
\cite{US} for details and \cite[Sec.~6]{BS-bio} for a recent treatment
in this context, with the transposed notation more common in
probability theory. We add two little examples to demonstrate that the
PBS can be used explicitly as well.

\section{Homogeneous case: Volterra equation and convergence}

Let $\pa (t; t^{}_0)$ denote the fundamental system of the homogeneous
part of \eqref{Cauchy}, which is also the flow of the linear
system. It satisfies the Volterra integral equation
\begin{equation} \label{FUSY}
   \pa (t; t^{}_{0}) \, = \, \one + \!\int_{t^{}_0}^{t} 
   \! A(\tau)\, \pa (\tau; t^{}_{0}) \dd\tau \ts ,
   \quad \text{with } \pa(t^{}_{0}; t^{}_{0}) = \one \ts .
\end{equation}
By means of a formal Picard iteration, this leads to
\begin{equation} \label{PBS}
    \pa (t; t^{}_0) \, = \, \one + \!\int_{t^{}_0}^{t} 
    \! A (\tau) \dd \tau \ts + \!\int_{t^{}_0}^{t} 
    \! A(\tau^{}_{1}) \int_{t^{}_0}^{\tau^{}_1} \! A(\tau^{}_{2})
      \dd \tau^{}_{2} \dd \tau^{}_{1} \, + \, \ldots   
      \; = \; \one + \sum_{n=1}^{\infty} \cI_{n} (t) \ts ,
\end{equation}
where
\[
   \cI_{n} (t) := \!\int_{t^{}_0}^{t} \! A(\tau^{}_{1})
   \int_{t^{}_0}^{\tau^{}_1} \!  A(\tau^{}_{2}) \, \cdots
   \int_{t^{}_0}^{\tau^{}_{n-1}}\! A(\tau^{}_{n}) \dd \tau^{}_{n} \cdots
   \dd \tau^{}_{2} \dd \tau^{}_{1} \ts .
\]
For convenience, we set $\cI^{}_{0} (t) = \one$. By construction, we
then have the recursion
\begin{equation} \label{recursion} 
   \cI_{n+1} (t) = \!
   \int_{t^{}_{0}}^{t} \! A(\tau) \, 
   \cI_{n} (\tau) \dd \tau
\end{equation}
for $n\in \NN_{0}$.  Clearly,
$\cI^{}_{n} (t^{}_{0}) =\delta^{}_{n,0} \one$, in line with
$\pa (t^{}_{0}; t^{}_{0})=\one$.  Eq.~\eqref{PBS} is known as the
Peano--Baker series (PBS) in control theory \cite[p.~598]{FH}, see
also \cite[Sec.~IV.5]{Gant}, or as Dyson's series expansion, compare
\cite[Rem.~6.2]{BS-bio}, in the context of the time-ordered
exponential in physics \cite{Hoek}. The PBS was recently also extended
to the setting of time scales \cite{Dac}.

Let us consider the individual terms of the PBS more closely.
\begin{lemma} \label{I-diffbar} When\/ $A$ is continuous on\/ $I$, the
  matrix functions\/ $\cI_{n}$ are continuously differentiable and
  satisfy\/ $\dot{\cI}_{n+1} (t) = A(t) \, \cI_{n} (t)$, for all\/
  $t\in I$ and\/ $n\in \NN_{0}$.
\end{lemma}

\begin{proof}
By definition \eqref{recursion}, each matrix function $\cI_{n}$ is
(component-wise) continuous, and hence so is the product
$A \cI_{n}$ for every $n$. Applying the fundamental theorem of
calculus separately to each component of the matrix equation
\eqref{recursion} proves the lemma.
\end{proof}

Let us show next (in modern terminology) that the PBS is nicely
convergent in our finite-dimensional setting; compare
\cite[Sec.~16.5]{I} or \cite[Sec.~2.11]{FDC} for a slighlty different
account, \cite{PDV} for further results in this direction, and
\cite[Ch.~4.3]{H} for background on matrix Taylor series. We begin by
establishing a link to the usual exponential series. In the
one-dimensional case, this is \cite[Cor.~1.3.1]{B}.

\begin{lemma} \label{dim-reduce} 
  Let\/ $A$ be continuous on the interval\/ $I$, with\/ $[A(t), A(s)]=0$ 
  for all\/ $t,s \in I$. Then, one has
\[
   \cI_{n} (t) = \frac{1}{n!} \Bigl(\, \int_{t^{}_0}^{t} 
   \! A(\tau) \dd\tau\Bigr)^n
\]   
   for all\/ $t\in I$ and\/ $n\in \NN^{}_{0}$. In particular, this applies
   when\/ $A$ is one-dimensional.
\end{lemma}

\begin{proof}
  The claim is obviously true for $n=0$ and $n=1$. We can now employ
  induction via the recursion \eqref{recursion}, assuming the validity
  for some $n\in \NN$. This gives
\[
   \cI_{n+1} (t) \, = \int_{t^{}_0}^{t} \! A(\tau^{}_{1})\, \cI_{n}
   (\tau^{}_{1}) \dd \tau^{}_{1} \, = \int_{t^{}_0}^{t} \!
   A(\tau^{}_{1})\, \frac{1}{n!}  \Bigl( \ts\int_{t^{}_0}^{\tau^{}_1} \!
   A(\tau) \dd\tau \Bigr)^n \dd\tau^{}_{1} \, ,
\]
where the integrand on the right hand side can be rewritten as
\[
   \frac{1}{n!} \ts A(\tau^{}_{1}) \Bigl( \, \int_{t^{}_0}^{\tau^{}_1} 
   \! A(\tau) \dd\tau \Bigr)^n 
   = \, \frac{1}{(n+1)!}\, \frac{\dd}{\dd \tau^{}_{1}}
   \Bigl( \, \int_{t^{}_0}^{\tau^{}_1} \! A(\tau) \dd\tau \Bigr)^{n+1}.
\]
This step employs the general chain rule, where the assumed
commutativity is used. Inserting this expression into the previous
formula completes the induction step by an application of the
fundamental theorem of calculus.
\end{proof}

\begin{remark}
When $[A(t), A(s)]=0$ on $I$, the PBS \eqref{PBS} reduces to the
well-known formula
\[
     \pa (t; t^{}_{0} ) \, = \sum_{n=0}^{\infty} \frac{1}{n!}\,
     \Bigl(\, \int_{t^{}_{0}}^{t} \! A(\tau) \dd \tau \Bigr)^n
     = \exp \Bigl(\, \int_{t^{}_{0}}^{t} \! A(\tau) \dd \tau \Bigr)
\]
as a consequence of Lemma~\ref{dim-reduce}; see \cite[Sec.~2.3]{KP}
for a detailed exposition of this case. This reduction is also
mentioned in most of the sources cited so far; compare
\cite[Cor.~1.3.2]{B}.
\end{remark}

\begin{remark}
  A closer look at the proof of Lemma~\ref{dim-reduce} shows that the
  condition $[A(t),A(s)]=0$ may be replaced by the slightly weaker
  assumption that, for all $t\in I$, the matrix $A(t)$ commutes with
  the integral $\int_{t^{}_0}^{t} \! A(\tau)\dd\tau$.

  Also, it should be mentioned that Lemma~\ref{dim-reduce} can
  alternatively be proved by direct calculations based on permutations
  of the integration variables followed by a suitable rearrangement to
  cover the integration region $[t^{}_{0},t]^{n}$.  When combined with
  an induction argument, it suffices to consider the permutations
  $(12\ldots n), \, (1)(23\ldots n), \ldots , \, (1)(2) \ldots
  (n\! - \!1,n)$. This approach is slightly more general. As it is
  also less transparent, we skip further details.
\end{remark}

Let now $\|.\|$ denote any norm on $\RR^d$ (or on $\CC^d$), and define
the compatible matrix norm by $\| A \| := \sup_{\|x\|=1} \| Ax\|$ as
usual. This implies the relations
$\| Ax \| \leqslant \| A\| \cdot \| x\|$ and, more importantly,
$\| AB \| \leqslant \| A\| \cdot \| B\|$.

\begin{theorem} \label{PBS-converges} 
  If\/ $\| A(t)\|$ is locally integrable on the interval\/ $I$, the
  series representation \eqref{PBS} of\/ $\pa (t; t^{}_{0})$ is
  compactly convergent on\/ $I$ in the chosen matrix norm.
\end{theorem}

\begin{proof}
  Let $J\subseteq I$ be compact, with $t_0\in J$. We show that the
  sequence of partial sums is Cauchy on $J$. So, let $m,n \in \NN$
  with $n > m$ and consider
\[
   \begin{split}
     \Bigl\| & \sum_{k=0}^{n} \cI_{k} (t) \ts - \sum_{k=0}^{m} \cI_{k}
     (t) \Bigr\| \; = \; \Bigl\| \sum_{k=m+1}^{n} \! \cI_{k}
     (t)\Bigr\| \; \leqslant \!
     \sum_{k=m+1}^{n} \bigl\| \ts \cI_{k} \bigr\| \\[1mm]
     & \leqslant \! \sum_{k=m+1}^{n} \int_{t^{}_{0}}^{t}
     \int_{t^{}_{0}}^{\tau^{}_{1}} \! \cdots
     \int_{t^{}_{0}}^{\tau^{}_{k-1}} \bigl\| A(\tau^{}_{1}) \cdots
     A(\tau^{}_{k}) \bigr\| \dd \tau^{}_{k} \cdots \dd \tau^{}_{2} \dd
     \tau^{}_{1} \ts .
   \end{split}  
\]
Since $\bigl\| A(\tau^{}_{1}) \cdots A(\tau^{}_{k}) \bigr\| \leqslant
\bigl\| A(\tau^{}_{1}) \bigr\| \cdots \bigl\| A(\tau^{}_{k}) \bigr\|$,
where all $\bigl\| A(\tau_i) \bigr\|$ are non-negative real numbers,
Lemma~\ref{dim-reduce} implies that the last sum is majorised by
\[
     \leqslant \sum_{k=m+1}^{n} \frac{1}{k!} \, \Bigl(\,
     \int_{t^{}_{0}}^{t} \bigl\| A(\tau) \bigr\| \dd \tau \Bigr)^k ,
\]
which is the corresponding Cauchy estimate for the Taylor series of
the ordinary exponential function on $\RR$, evaluated at $
\int_{t^{}_{0}}^{t} \bigl\| A(\tau) \bigr\| \dd \tau$, which exists
for all $t\in J$ by assumption. Since this series converges compactly,
our claim follows.
\end{proof}

\section{Solution of the inhomogeneous problem}

It is now obvious that $\pa (t; t^{}_{0})$ solves the homogeneous
initial value problem \eqref{FUSY}. This follows from a term-wise
application of Lemma~\ref{I-diffbar} to the PBS, which is fully
justified by Theorem~\ref{PBS-converges}. The determinant of $\pa$,
which is a Wronskian and thus satisfies Liouville's theorem, reads
\begin{equation} \label{Wronski} 
  \det \bigl( \pa (t; t^{}_{0}) \bigr)
  \, = \: \det \bigl( \pa (t^{}_{0}; t^{}_{0}) \bigr) \nts\nts \cdot \exp \Bigl( \,
  \int_{t^{}_{0}}^{t} \mathrm{tr} \bigl( A(\tau) \bigr) \dd \tau
  \Bigr) \, = \: \exp \Bigl( \, \int_{t^{}_{0}}^{t} \mathrm{tr} \bigl(
  A(\tau) \bigr) \dd \tau \Bigr) ,
\end{equation}  
which never vanishes; see \cite[Prop.~11.4]{A} or \cite[Thm.~2.23]{KP}
for details.  This means that $\pa$ has full rank and thus indeed
constitutes a fundamental system of the homogeneous linear system.
  
Since $\pa (t; t^{}_{0})$ is the unique solution of \eqref{FUSY}, the
flow property implies the relation
\begin{equation} \label{flow-prop}
   \pa (t; s) \, \pa (s; t^{}_{0}) \, = \ts \, \pa (t; t^{}_{0}) \ts ,
\end{equation}  
which, due to $\pa (t; t)=\one$, also implies $\bigl( \pa (t; s)
\bigr)^{-1} = \pa (s; t)$.  With the usual `variation of constants'
calculation, compare \cite[Thm.~11.13]{A}, one can now easily derive
the following result.

\begin{theorem} \label{gen-Cauchy} 
  Let\/ $I$ be an interval, with\/ $t^{}_{0} \in I$. Suppose\/ $A$ is a
  continuous matrix function on\/ $I$, and\/ $b$ is also continuous on\/
  $I$. Then, the Cauchy problem \eqref{Cauchy} has the unique solution
\[
    x(t) \, = \,  \pa (t; t^{}_{0}) \Bigl( x^{}_{0} + \int_{t^{}_{0}}^{t}
    \pa (t^{}_{0}; \tau)\, b(\tau) \dd \tau \Bigr)
\]  
with\/ $\pa$ given by the PBS \eqref{PBS}.  When\/ $[A(t),A(s)]=0$ 
for all\/ $t,s\in I$, the formula simplifies to the explicit expression
\eqref{commuting-case} with ordinary exponentials. \qed
\end{theorem}

\section{Examples}

Let us demonstrate the explicit applicability of the PBS with two simple
examples. Both can also be solved by other means, but are perhaps
still instructive.

Consider the matrix family $A(t) = \left( \begin{smallmatrix} 1 & t \\
    0 & a \end{smallmatrix} \right)$ with $t\in I$, $t^{}_{0} = 0$,
and fixed parameter $a$.  The matrices commute for $a=1$, but not
otherwise. One finds
\[
    \cI_{n} (t) \, = \, \begin{pmatrix}
    \frac{t^{n}}{n!} & \frac{t^{n+1}}{(n+1)!}\, \alpha^{}_{n} \\
    0 &  \frac{\vphantom{\hat{\hat{T}}}(a t)^{n}}{n!} \end{pmatrix}
    \quad \text{with} \quad \alpha^{}_{n} =
    \sum_{\ell=1}^{n} \ell \ts a^{\ell-1} ,
\]
so that the PBS gives
\[
    \pa (t; 0) \, = \, \begin{pmatrix}
    e^t & f(t) \\ 0 & e^{at\vphantom{\hat{t}}} \end{pmatrix}
    \quad \text{with} \quad f(t) =
    \frac{e^{t} - e^{at} - (1-a)\ts t \ts e^{at}}{(1-a)^{2}} \, .
\]
Note that $f(t)$ simplifies to $\frac{1}{2} t^{2} \ts e^{t}$ for
$a=1$, in line with the then simpler ODE system.  For general $t_{0}$,
one determine $\pa (t;t_{0})$ from a similar calculation.  The PBS
differs both from the matrix exponential and from the known Mathias formula
for upper-triangular matrices; compare \cite[Thm.~3.6]{H}.

As a second example, consider the ODE system
\[
   \begin{pmatrix} \dot{x} \\ \dot{y} \end{pmatrix} \, = \,
   \begin{pmatrix} 0 & t \\ a & 0\end{pmatrix}
   \begin{pmatrix} x \\ y \end{pmatrix} ,
\]
which leads to the Airy function via the 2nd order ODE $\ddot{y} = a\ts  t
\ts y$; see \cite[Ch.~10.4]{AS} for details. With $\alpha := a^{1/3}$,
the PBS leads to the formula
\[
    \pa (t;0) \, = \,
    \begin{pmatrix}  \dot{g} (\alpha t) & \frac{1}{\alpha^{2} }
    \dot{f} (\alpha t) \\
    \alpha^{2} g(\alpha t) & f(\alpha t)    \end{pmatrix}
\]
with
\[
     f(z) \, = \sum_{k=0}^{\infty} 3^{k} (\tfrac{1}{3})^{}_{k} 
     \frac{z^{3k}}{(3k)!}
     \quad \text{and} \quad
     g(z) \, = \sum_{k=0}^{\infty} 3^{k} (\tfrac{2}{3})^{}_{k} 
     \frac{z^{3k+1}}{(3k+1)!} 
\]
from \cite[Eq.~10.4.3]{AS}.  Here, $(\beta)^{}_{0} = 1$ and
$(\beta)^{}_{k} = \beta \ts (\beta+1) \cdots (\beta + k - 1)$ for
$k\in\NN$ as usual; compare \cite[Eq.~6.1.22]{AS}. Note that
\[
    \pa (t; 0) \; \xrightarrow{\, a \to 0 \,} \; \begin{pmatrix}
    1 & \frac{t^{2}}{2} \\ 0 & 1 \end{pmatrix},
\]
in line with the trivially solvable system for $a=0$.

\section*{Acknowledgements}
It is a pleasure to thank Peter Jarvis for discussions, and
Wolf-J\"urgen Beyn and Nikolai Dolbilin for useful suggestions.

\end{document}